\documentclass[12pt]{article}
\usepackage{amsmath,amssymb,latexsym}
\usepackage{array,graphicx,hyperref}
\usepackage{hyperref,stmaryrd,xcolor,psfrag} 
\usepackage{pgf,tikz,pgfkeys,pgfplots,float}
\usepackage[latin1]{inputenc}

\newtheorem{theo}{Theorem}
\newtheorem{prop}{Proposition}
\newtheorem{rk}{Remark}

\newcommand{\proof}{\medskip\goodbreak\noindent {\textbf{Proof.} }}
\def\endproof{\hfill{}$\square$} 

\newcommand\RR{\mathbb{R}}

\providecommand{\abs}[1]{\lvert#1\rvert}
\providecommand{\norm}[1]{\lVert#1\rVert}
\newcommand\ds{\displaystyle}

\begin{document}

\title{Two approaches for the stabilization of  nonlinear KdV equation with boundary time-delay feedback.}

\author{Lucie Baudouin\footnotemark[1]  
\and Emmanuelle Cr\'epeau\footnotemark[2]
\and Julie Valein\footnotemark[3]
}

\maketitle
\footnotetext[1]{LAAS-CNRS, Universit\'e de Toulouse, CNRS, UPS, Toulouse, France. {\tt\small baudouin@laas.fr}}
\footnotetext[2]{LMV, UVSQ, CNRS,  Universit\'e Paris-Saclay, 78035 Versailles, France. {\tt\small emmanuelle.crepeau@uvsq.fr}}
 \footnotetext[3]{Institut Elie Cartan de Lorraine, Universit\'e de Lorraine \& Inria (Project-Team SPHINX), BP 70 239, F-54506 Vand{\oe}uvre-les-Nancy Cedex, France. {\tt\small julie.valein@univ-lorraine.fr    }}
\footnotetext{This work was partially supported by  MathAmsud project, ICoPS and the ANR project SCIDiS contract number 15-CE23-0014.}

\maketitle

\begin{abstract}    
This article concerns the nonlinear Korteweg-de Vries equation with boundary time-delay feedback. Under appropriate assumption on the  coefficients of the feedbacks (delayed or not), we first prove that this nonlinear infinite dimensional system is well-posed for small initial data. The main results of our study are two theorems stating the exponential stability of the nonlinear time delay system. Two different methods are employed: a Lyapunov functional approach (allowing to have an estimation on the decay rate, but with a restrictive assumption on the length of the spatial domain of the KdV equation) and an observability inequality approach, with a contradiction argument (for any non critical lengths but without estimation on the decay rate). Some numerical simulations are given to illustrate the results.
\end{abstract}

\section{Introduction and main results}

In 1834, John Scott Russell observed for the first time in a Scottish canal a solitary wave, also called soliton, which propagates without deformation in a nonlinear and dispersive medium. We refer to the very good introduction of \cite{RosierZhang09} for the literary description of this phenomenon. 
In 1895, Diederik Korteweg and Gustav de Vries derived the nonlinear dispersive partial differential equation which is now known as the Korteweg-de Vries (KdV) equation and models the propagation of a long wave in water of relatively shallow depth. 
It seems that this equation was first introduced by Joseph Boussinesq in 1877, \cite{Boussinesq}.
The domains of applications of this equation are various: collision of hydromagnetic waves, ion acoustic waves in a plasma, acoustic waves on a crystal lattice or even subparts of the cardiovascular system... 
We refer for instance to the book \cite{book:Whitham} for a physical deduction of the KdV equation.\\

Adding a delayed term in the boundary stabilization of this equation is a way to  take into account the reality of any device placed to implement a boundary feedback control. A domain of application could be for instance the study of pulsatile flow in blood vessels that can be modelled by the KdV equation (see \cite{Sorine}),  but we will work here specifically on the mathematical technicalities of the proof of stabilization results stated below.
The challenge on the specific topic of our contribution, beyond the difficulty of dealing with a nonlinear equation, is to prove that under appropriate conditions, a delay in the boundary feedback of this equation will not destabilize the system \cite{Datko_1988}.

The first work concerning the exponential stabilization of the KdV equation (without delay) on a bounded domain is \cite{Zhang_1994}, where the length of the spatial domain is $L=1$ and holds under an appropriate assumption on the weight of the feedback. 
As expected when dealing with the KdV equation, the length $L$ of the domain where the equation is set plays a role in the ability of controlling (\cite{Cerpa_2014}, \cite{Coron_2004}, \cite{Rosier_1997}) or stabilizing (\cite{BK-TAC00}, \cite{Perla_2002}, \cite{RosierZhang09}) the solution of the equation. Indeed, it is well-known that if $L=2\pi$, there exists a solution ($y(x,t)=1-\cos x$) of the linearized system around~$0$ which has a constant energy. More generally,
defining the set of critical lengths
$$\mathcal N=\left\{2\pi\sqrt{\frac{k^2+kl+l^2}{3}}, \, k,l\in\mathbb N^*\right\},$$
one can recall that the linearized equation around $0$  is exactly controllable if and only if $L\notin \mathcal N$ (see \cite{Rosier_1997}) and the local exact controllability result holds for the nonlinear KdV equation (using a fixed point argument) if $L\notin \mathcal N$. Further results show that the nonlinear KdV equation is in fact locally exactly controllable for all critical lengths contrary to the linear KdV equation (see \cite{Coron_2004},~\cite{Cerpa2007},~\cite{CerpaCrepeau_2009}). 

Nevertheless, concerning the stabilization topic, in the case of non critical length, it is not necessary (see \cite{Perla_2002}) to introduce a feedback law as in \cite{Zhang_1994} to have the local exponential stability of the nonlinear KdV equation. Moreover, it is proved in \cite{Perla_2002} and \cite{Pazoto_2005} that for any critical length, adding a localized damping in the nonlinear KdV equation allows to have a local exponential stability result, and even a semi-global stability by working directly with the nonlinear system. Recently, in \cite{chu2015asymptotic} and \cite{tang2016asymptotic} some results of asymptotic stability for the nonlinear KdV equation for the first critical length ($2\pi$) and the second one ($2\pi\sqrt{\frac{7}{3}}$) have been proven without any feedback law.

Even if the exponential stability holds, a related interesting question is the rapid stabilization, or how to construct a feedback law which stabilizes the system at a prescribed decay rate. In this context, we should mention \cite{KS-SCL08} that deals with the stabilization of a linear KdV-like equation with the backstepping method, \cite{CerpaCoron_2013} for the local rapid stabilization for the nonlinear KdV equation from the left Dirichlet boundary condition by the same method and \cite{CoronLu_2014} for the local rapid stabilization for a KdV equation with a Neumann boundary control on the right by an integral transform for non critical lengths. Finally we also refer to the recent paper \cite{Marx_2017} for  the global stabilization of a nonlinear KdV equation with a saturating distributed control.\\ 

{\bf Notations:}
 $L^2(a,b)$ represents the space of square integrable functions over the interval $(a,b)$ with values in $\mathbb R$ and the partial derivatives in time and space of a function $y$ are denoted $y_t$ and $y_x$. The Hilbert space $H^1(a,b)$ (resp. $H^3(a,b)$) is the set of all functions $y\in L^2(a,b)$ such that $y_x\in L^2(a,b)$ (resp. $y_x$, $y_{xx}$ and $y_{xxx} \in L^2(a,b)$). We will also use (and recall) the following functional spaces :
 $\mathcal H= L^2(0,L)\times L^2(-h,0)$, $\mathcal B= C([0,T],L^2(0,L))\cap L^2(0,T,H^1(0,L))$ and   $H=L^2(0,L)\times L^2(0,1)$.\\

The main goal  of this paper is to study the stabilization of the following nonlinear KdV equation with a boundary feedback delayed term
\begin{equation}\label{syst_nonlinear}
\hspace{-0.3cm}\left\{\hspace{-0.1cm}\begin{array}{l}
y_t(x,t)+y_{xxx}(x,t)+y_x(x,t)+y(x,t)y_x(x,t)=0, \\
 \hfill{~} x\in(0,L),\,t>0,\\
y(0,t)=y(L,t)=0,  \hfill{~}  t>0,\\
y_x(L,t)=\alpha y_x(0,t)+\beta y_x(0,t-h),  \hfill{~}  t>0,\\
y_x(0,t)=z_0(t),  \hfill{~} t\in(-h,0),\\
y(x,0)=y_0(x),  \hfill{~}  x\in(0,L),
\end{array}\right.
\end{equation}
where $h>0$ is the delay, $L>0$ is the length of the spacial domain, $\alpha$ and $\beta\neq 0$ are real constant parameters and $y(x,t)$ is the amplitude of the water wave at position $x$ at time~$t$. The initial data $y_0$ is supposed to belong to $L^2(0,L)$ and the delayed left lateral boundary Neumann data $z_0$ belongs to $L^2(-h,0)$.

We define the Hilbert space of the initial and boundary data $\mathcal H := L^2(0,L)\times L^2(-h,0)$, 
endowed with the norm defined for all $(y,z)\in \mathcal H $ by
$$\left\|(y,z)\right\|_{\mathcal H}^2 = \int_0^L y^2(x) dx + \abs{\beta}h\int_{-h}^0 z^2(s) ds.
$$
In the case without delay (i.e. $\beta=0$) it is well-known (see for instance \cite{Zhang_1994}) that for every $T>0$, $L>0$ and $y_0\in L^2(0,L)$, the system \eqref{syst_nonlinear} is locally well-posed in 
$C([0,T],L^2(0,L))\cap L^2(0,T,H^1(0,L)) := \mathcal B$. 
We will give in Section~\ref{sect:wellpo} the proof of well-posedness for the case with delayed boundary condition (i.e. $\beta \neq 0$).\\

Before stating the two main results of this article, let us precise what is at stake when working with the nonlinear KdV equation. Knowing that critical values $\mathcal N$ of the length $L$ of the spatial domain are precluding the system from controllability and stabilizability as soon as the study involves the linearized equation around $0$ (see \cite{Rosier_1997}), we are expecting the best  possible stabilization result for $L$ no larger than the first critical length $2\pi$. 
Let us give the following definition of the energy of system~\eqref{syst_nonlinear}, chosen because it corresponds to the norm of 
$(y(\cdot, t), y_x(0,t-h~\cdot))$ on $\mathcal H$:
\begin{equation}\label{def:E}
E(t)=\int_0^L y^2(x,t) dx + \abs{\beta}h\int_0^1 y_x^2(0,t-h\rho) d\rho.
\end{equation}
One should know that this is a classical choice when considering boundary delayed terms, as in \cite{Nicaise_2006} and \cite{NVF-DCDS09} for the heat and wave partial differential equations.
Moreover, we will assume, along the whole article, that the coefficients $\alpha$ and $\beta$ comply to  the following limitation:
\begin{equation}\label{hyp:alpha_beta}
\left|\alpha\right|+\left|\beta\right|<1.
\end{equation}
This is necessary even for the existence of solutions and we can refer to the case without delay ($\beta = 0$, $|\alpha| <1 $) that can be read in \cite{Zhang_1994}.
Note also that, in \cite{Nicaise_2006}, dealing with the wave equation, there are some restrictions about the positive coefficients of the terms with or without delay. Actually, it is the case for hyperbolic and parabolic partial differential equations in \cite{NV2010}, \cite{NVF-DCDS09} and even for the Schr\"odinger equation (which is a dispersive equation, just like KdV) in \cite{NR-PM11}. In these papers, the authors assume that the coefficient of the term with delay is smaller than the coefficient of the term without delay, i.e. with our notations here : $0\leq\beta<\alpha$.
This kind of assumption is necessary in these cases and if they are not satisfied, it can be shown that instabilities may appear (see for instance \cite{Datko_1988}, \cite{Datko_1986} with $\alpha=0$, or \cite{Nicaise_2006} in the more general case for the wave equation).
For the KdV equation we do not have this kind of assumption. We can take $\alpha=0$ (and then $\left|\beta\right|<1$) or even $\alpha=\beta=0$, adapting of course the inner product and the proofs of the main results. The main goal of this paper is to show that a delay does not destabilize the system, contrary to many other delayed systems (see \cite{Datko_1988}, \cite{Datko_1986}, \cite{Logemann_1996}, \cite{Rebarber_1999}, \cite{GY-TAC10}).\\

Our first main result is obtained for a restricted assumption on the length $L$ but yields exponential stability of the solution of system~\eqref{syst_nonlinear} with an estimation of the decay rate stated below.

\begin{theo}\label{thm:expostab_nonlinear}
{ \it Assume that $\alpha$ and $\beta \neq 0$ satisfy \eqref{hyp:alpha_beta} and assume that the length $L$ fulfills
\begin{equation}\label{hyp:L}
L<\pi\sqrt{3}.
\end{equation} 
Then, there exist $r>0$, $\mu_1>0$ and $\mu_2\in(0,1)$ sufficiently small, such that for every $(y_0,z_0)\in \mathcal H $ satisfying
$$\left\|(y_0,z_0)\right\|_{\mathcal H}\leq r,$$ 
the energy of system~\eqref{syst_nonlinear}, denoted $E$ and defined by \eqref{def:E}, decays exponentially. More precisely, there exist two positive constants $\gamma$ and $\kappa$ such that 
$$
E(t)\leq \kappa E(0)e^{-2\gamma t},\qquad t>0,
$$
where 
\begin{equation}\label{eq:gamma}
\hspace{-0.2cm}\gamma\leq\min\left\{\frac{(9\pi^2-3L^2-2L^{3/2}r\pi^2)\mu_1}{6L^2(1+L\mu_1)},\frac{\mu_2}{2(\mu_2+\abs{\beta})h}\right\}.
\end{equation}
}
\end{theo}
This theorem will be proved in a constructive manner, allowing an estimation of the decay rate~$\gamma$. Up to our knowledge, such a quantitative estimation is new in the literature of non-linear KdV stabilization. The proof  uses an appropriate Lyapunov functional build with coefficients $\mu_1$ and $\mu_2$ and detailed in Section~\ref{LyapSec}.
Moreover, note that when the delay $h$ becomes larger, then the decay rate~$\gamma$ is smaller.

\begin{rk}
The coefficients $\mu_1$ and  $\mu_2$ depend on the Lyapunov functional we will use in the proof of the stability result. Nevertheless, one can have an estimation of both $r$, $\mu_1$ and $\mu_2$ in Remark \ref{rk:mubounded} below.
\end{rk}

On the other hand, our second main result is obtained simply for non critical lengths and gives generic exponential stability of the solution of system \eqref{syst_nonlinear}.
\begin{theo}\label{thm:expostab_nonlinear_LN}
{ \it Assume that the length $L>0$ satisfies $L \notin\mathcal N$ and that $\alpha$ and $\beta$ satisfy \eqref{hyp:alpha_beta}.
Then, there exists $r>0$ such that for every $(y_0,z_0)\in \mathcal H$ satisfying
$$\left\|(y_0,z_0)\right\|_{\mathcal H}\leq r,$$ 
the energy of system~\eqref{syst_nonlinear}, denoted $E$ and defined by \eqref{def:E}, decays exponentially. More precisely, there exist two positive constants $\nu$ and $\kappa$ such that 
$$
E(t)\leq \kappa E(0)e^{-\nu t},\qquad t>0.
$$
}
\end{theo}
The proof of this theorem relies on an observability inequality and the use of a contradiction argument. Thus, the value of the decay rate can not be estimated precisely in this approach.\\

These two results of Theorem~\ref{thm:expostab_nonlinear} and \ref{thm:expostab_nonlinear_LN} both have a specific interest and a dedicated methodological approach that are worth being presented one after another. The next section is devoted to the necessary preliminary step dealing with the well-posedness and regularity of the solutions of our specific system coupling  the KdV equation and a delayed boundary feedback. Section~\ref{LyapSec} will develop the proof of a first quantified exponential stabilization result stated in Theorem~\ref{thm:expostab_nonlinear} while the proof of our second  stabilization result, stated in Theorem~\ref{thm:expostab_nonlinear_LN}, will be detailed in Section~\ref{ObsSec}. When necessary, a first step concerning the linearized KdV equation will be given. Finally, Section~\ref{NumSec} will detail a numerical simulation meant to illustrate our work.

\section{Well-posedness and regularity results}\label{sect:wellpo}

\subsection{Study of the linear equation}
We begin by proving the well-posedness of the KdV equation linearized around 0, that writes
\begin{equation}\label{syst_linear}
\hspace{-0.3cm} \left\{\hspace{-0.1cm}\begin{array}{ll}
y_t(x,t)+y_{xxx}(x,t)+y_x(x,t)=0,  \hfill{~} \! x\! \in \!(0,L), t>0,\\
y(0,t)=y(L,t)=0,  \hfill{~}  t>0,\\
y_x(L,t)=\alpha y_x(0,t)+\beta y_x(0,t-h),  \hfill{~}  t>0,\\
y_x(0,t)=z_0(t),  \hfill{~}  t\in(-h,0),\\
y(x,0)=y_0(x),  \hfill{~}  x\in(0,L).
\end{array}\right.
\end{equation}
Following Nicaise and Pignotti \cite{Nicaise_2006}, we set $z(\rho,t)=y_x(0,t-\rho h)$ for any $\rho\in(0,1)$ and $t>0$. Then $z$ satisfies the transport equation
\begin{equation}\label{syst_z}
\left\{\begin{array}{ll}
h z_t(\rho,t)+z_{\rho}(\rho,t)=0, & \rho\in(0,1),\,t>0,\\
z(0,t)=y_x(0,t), & t>0,\\
z(\rho,0)=z_0(-\rho h), & \rho\in(0,1).
\end{array}\right.
\end{equation}
We introduce  the Hilbert space $L^2(0,L)\times L^2(0,1) := H$ equipped with the inner product
$$\left\langle\left(\begin{array}{c}y\\ z\end{array}\right), \left(\begin{array}{c}\tilde y\\ \tilde z\end{array}\right) \right\rangle
= \int_0^L y\tilde y \,dx + \abs{\beta}h \int_0^1 z\tilde z \,d\rho,$$ 
for any $(y, z), (\tilde y, \tilde z)\in H$. 
This new inner product is clearly equivalent to the usual inner product on $H$ and we denote by $\left\|\cdot\right\|_H$ the associated norm.

We then rewrite \eqref{syst_linear} as a first order system:
\begin{equation}\label{EqA}
\left\{\begin{array}{ll}
U_t(t)=\mathcal A U(t), & t>0,\\
U(0)=U_0\in H,&
\end{array}\right.
\end{equation}
 where 
$U=\left(\begin{array}{c}y\\ z\end{array}\right)$, 
$U_0=\left(\begin{array}{c}y_0\\ z_0(- h~\cdot)\end{array}\right)$,
and where the operator is defined by
$$\mathcal A =\left(\begin{array}{cc}-\partial_{xxx}-\partial_x& 0\\ 0&-\frac{1}{h}\,\partial_{\rho}\end{array}\right),$$
with domain
$$
\mathcal D(\mathcal A)=\Big\{(y,z)\in H^3(0,L)\times H^1(0,1)\left| \,y(0)=y(L)=0\right.,
 z(0)=y_x(0), y_x(L)=\alpha y_x(0)+\beta z(1)\Big\}.
$$

\begin{theo}\label{thm:wellposed}
{ \it Assume that $\alpha$ and $\beta$ satisfy \eqref{hyp:alpha_beta} and that  $U_0\in H$. Then there exists a unique mild solution $U\in C([0,+\infty),H)$ for system \eqref{EqA}. Moreover if $U_0\in\mathcal D(\mathcal A)$, then the solution is classical and satisfies 
$$U\in C([0,+\infty),D(\mathcal A))\cap C^1([0,+\infty),H).$$
}
\end{theo}

\proof
We first prove that the operator $\mathcal A$ is dissipative.
Let $U=(y,z)\in \mathcal D(\mathcal A)$.
Then we have
$$
	\begin{aligned}
	&\left\langle \mathcal A U, U\right\rangle =  \ds -\int_0^L y_{xxx} y \,dx -\int_0^L y_x y \,dx 
- \abs{\beta} \int_0^1 z_\rho z \,d\rho\\
& = \frac12 [y_x^2]_0^L - \frac{\abs{\beta}}{2} [z^2]_0^1\\
& =   \frac{1}{2} \left(\alpha y_x(0)+\beta z(1)\right)^2- \frac{1}{2} y_x^2(0)- \frac{\abs{\beta}}{2} z^2(1) + \frac{\abs{\beta}}{2} y_x^2(0)\\
& =   \frac12 \left( M \xi,\xi\right),
	\end{aligned}
$$
where $(.,.)$ is the usual scalar product in $\RR^2$,
\begin{equation}\label{def:M}
\xi=\begin{bmatrix}y_x(0)\\ z(1)\end{bmatrix} \hbox{ and } 
M=\begin{bmatrix}\alpha^2-1+\abs{\beta} & \alpha\beta \\ \alpha\beta & \beta^2-\abs{\beta}\end{bmatrix}.
\end{equation}
One can check that $M$ is definite negative: the trace of $M$ satisfies 
${\rm tr}  M=\alpha^2+\beta^2-1<0$ 
if and only if  $\alpha^2+\beta^2<1$,
and the determinant of $ M$ gives
$\det  M = \abs{\beta}\left((\abs{\beta}-1)^2-\alpha^2\right)$,
so that $ M$ is definite negative iff ($\alpha^2+\beta^2<1$ and $(\abs{\beta}-1)^2-\abs{\alpha}^2>0$), which is equivalent to hypothesis \eqref{hyp:alpha_beta}.

Secondly we show that the adjoint of $\mathcal A$, denoted by $\mathcal A^*$, is also dissipative.
It is not difficult to prove that the adjoint is defined by
$\mathcal A^* =\left(\begin{array}{cc}\partial_{xxx}+\partial_x& 0\\ 0&\frac{1}{h}\,\partial_{\rho}\end{array}\right),$
with domain 
\begin{multline*}
\mathcal D(\mathcal A^*)=\Big\{(\tilde y,\tilde z)\in H^3(0,L)\times H^1(0,1)\big|\tilde y(0)=\tilde y(L)=0,\\ 
\tilde z(1)=\frac{\beta}{\abs{\beta}}\, \tilde y_x(L) \hbox{ and } \tilde y_x(0)=\alpha \tilde y_x(L)+\abs{\beta} \tilde z(0)\Big\}.
\end{multline*}
Then for all $\tilde U=(\tilde y,\tilde z)\in \mathcal D(\mathcal A^*)$,
$$
\left\langle \mathcal A^* \tilde U, \tilde U\right\rangle  =  \ds \int_0^L \tilde y_{xxx} \tilde y \,dx +\int_0^L\tilde y_{x} \tilde y \,dx 
+ \abs{\beta} \int_0^1 \tilde z_\rho \tilde z \,d\rho
 = ~ \ds -\frac12 [\tilde y_x^2]_0^L + \frac{\abs{\beta}}{2} [\tilde z^2]_0^1
 =  \ds \frac12 \left(\tilde M \tilde \xi,\tilde \xi\right),
$$
where
$\tilde \xi=\begin{bmatrix}\tilde y_x(L)\\ \tilde z(0)\end{bmatrix}$  and 
$\tilde M =\begin{bmatrix}\alpha^2+\abs{\beta}-1 & \alpha \abs{\beta} \\ \alpha \abs{\beta} & \beta^2-\abs{\beta}\end{bmatrix}$.
Since ${\rm tr} \tilde M ={\rm tr}  M$ and $\det \tilde M=\det M$, we deduce that under hypothesis \eqref{hyp:alpha_beta}, $\tilde M$ is definite negative.

Finally, since $\mathcal A$ is a densely defined closed linear operator, and both $\mathcal A$ and $\mathcal A^*$ are dissipative, then $\mathcal A$ is the infinitesimal generator of a $C_0$ semigroup of contractions on $H$, which finishes the proof. 
\endproof\\

We denote by $\left\{S(t),t\geq0\right\}$ the semigroup of contractions associated with $\mathcal A$. In the following the real $C$ designs a positive constant that can depend on $T,\beta, h$. 
Let us now detail a few a priori estimates and regularity estimates of the solutions of system~\eqref{syst_linear}.

\begin{prop}\label{regularity}
{ \it Assume that \eqref{hyp:alpha_beta} is satisfied. Then, the map 
\begin{eqnarray}
(y_0,z_0(-h~\cdot)) & \mapsto & S(\cdot)(y_0,z_0(-h~\cdot)) \label{ContMap}
\end{eqnarray}
is continuous from $H$ to $\mathcal B \times C([0,T],L^2(0,1))$,
and for $(y_0,z_0(-h~\cdot))\in H$, one has \\
$(y_x(0,.),z(1,.))\in (L^2(0,T))^2$ and the following estimates
\begin{eqnarray}
\norm{y_x(0,\cdot)}^2_{L^2(0,T)}+\norm{z(1,\cdot)}^2_{L^2(0,T)} 
\leq C\left(\norm{y_0}^2_{L^2(0,L)}+\norm{z_0(-h~\cdot)}^2_{L^2(0,1)}\right), \label{y_x(0)+z(1)}\\
\norm{y_0}^2_{L^2(0,L)}\leq \frac{1}{T}\norm{y}^2_{L^2(0,T,L^2(0,L))}+\norm{y_x(0,\cdot)}^2_{L^2(0,T)}\label{y_0},\\
\norm{z_0(-h\,\cdot)}^2_{L^2(0,1)}\leq \norm{z(\cdot,T)}^2_{L^2(0,1)} \!+ \!\frac{1}{h}\norm{z(1,\cdot)}^2_{L^2(0,T)}.\label{z_0}
\end{eqnarray}
}
\end{prop}

\proof

\noindent $\bullet$ First of all, for any $(y_0,z_0(-h~\cdot))\in H$, Theorem \ref{thm:wellposed} brings $S(.)(y_0,z_0(-h~\cdot))=(y,z)\in C([0,T],H)$ and as the operator $\cal A$ generates a $C_0$ semi-group of contractions we get for all $t\in [0,T]$,
\begin{multline}\label{norm_H}
\norm{y(t)}^2_{L^2(0,L)}+\abs{\beta}h\norm{z(t)}^2_{L^2(0,1)}\\
\leq \norm{y_0}^2_{L^2(0,L)}+\abs{\beta}h\norm{z_0(-h~\cdot)}^2_{L^2(0,1)}.
\end{multline}

Let $p\in C^{\infty}([0,1]\times[0,T])$, $q\in C^{\infty}([0,L]\times[0,T])$ and $(y,z)\in \cal D(\cal A)$. Then multiplying \eqref{syst_z} by $pz$ and \eqref{syst_linear} by $qy$, and using some integrations by parts we get
\begin{multline}\label{eq_p_z}
\int_0^1 \left(p(\rho,T)z^2(\rho,T)-p(\rho,0)z_0^2(-\rho h)\right)d\rho
-\frac{1}{h}\int_0^T\!\!\int_0^1(hp_t+p_\rho)z^2d\rho dt\\
+\frac{1}{h}\int_0^T \left(p(1,t)z^2(1,t)-p(0,t)y_x^2(0,t)\right)dt=0
\end{multline}
\begin{multline}\label{eq_q_y}
\int_0^L \left(q(x,T)y^2(x,T)-q(x,0)y_0^2(x)\right)dx
-\int_0^T\!\!\int_0^L(q_t+q_x+q_{xxx})y^2dxdt+3\int_0^T\!\!\int_0^Lq_xy_x^2dxdt\\
-\int_0^Tq(L,t)(\alpha y_x(0,t)+\beta z(1,t))^2dt
+\int_0^T q(0,t)y_x^2(0,t)dt=0.
\end{multline}

\noindent $\bullet$ Let us first choose $p(\rho,t)\equiv\rho$ in \eqref{eq_p_z}. Then we obtain 
$$
\int_0^1 \rho\left(z^2(\rho,T)-z_0^2(-\rho h)\right)d\rho
-\frac{1}{h}\int_0^T\int_0^1z^2d\rho dt
+\frac{1}{h}\int_0^T z^2(1,t)dt=0
$$
and thanks to \eqref{norm_H} we get
\begin{equation}\label{ineq_z}
\hspace{-0.2cm}\|z(1,.)\|^2_{L^2(0,T)}\leq C\big(\norm{y_0}^2_{L^2(0,L)}+\norm{z_0(-h~\cdot)}^2_{L^2(0,1)}\big).
\end{equation}

Secondly, if we choose $q(x,t)\equiv1$ in \eqref{eq_q_y}, then we get,
$$
\int_0^L \left(y^2(x,T)-y_0^2(x)\right)dx
-\int_0^T(\alpha y_x(0,t)+\beta z(1,t))^2dt + \int_0^T y_x^2(0,t)dt = 0,
$$
which implies
$$
\ds \int_0^T y_x^2(0,t)dt\leq \ds \int_0^T(\alpha y_x(0,t)+\beta z(1,t))^2dt + \|y_0\|^2_{L^2(0,L)}.
$$
Therefore, since 
$$(\alpha y_x(0,t)+\beta z(1,t))^2\leq (\alpha^2+\beta^2)\left(y_x^2(0,t)+z^2(1,t)\right),$$
 we obtain
$$
\int_0^T(1-(\alpha^2+\beta^2))y_x^2(0,t)dt
\leq \int_0^T(\beta^2+\alpha^2)z^2(1,t)dt + \|y_0\|^2_{L^2(0,L)}
$$
and using the previous estimate of $\|z(1,.)\|_{L^2(0,T)}$ and hypothesis \eqref{hyp:alpha_beta}, we get
$$
\|y_x(0,.)\|^2_{L^2(0,T)}\leq C\left(\norm{y_0}^2_{L^2(0,L)}+\norm{z_0(-h~\cdot)}^2_{L^2(0,1)}\right)
$$ 
that concludes the proof of \eqref{y_x(0)+z(1)}.

 \noindent $\bullet$ Taking  now  $q(x,t)\equiv x$ in \eqref{eq_q_y}, we can write 
 \begin{multline*}
 \int_0^L x\left(y^2(x,T)-y_0^2(x)\right)dx-\int_0^T\int_0^Ly^2dxdt\\
 +3\int_0^T\int_0^Ly_x^2dxdt
-\int_0^T L(\alpha y_x(0,t)+\beta z(1,t))^2dt=0
 \end{multline*}
and using \eqref{norm_H} and \eqref{y_x(0)+z(1)} we obtain 
$$
	\|y_x\|^2_{L^2(0,T;L^2(0,L))}\leq C\left(\norm{y_0}^2_{L^2(0,L)}+\norm{z_0(-h~\cdot)}^2_{L^2(0,1)}\right)
$$
that brings, together with \eqref{norm_H}, the continuity of the map \eqref{ContMap}.
  
  \noindent $\bullet$ Choosing $q(x,t)\equiv T-t$ in \eqref{eq_q_y} yields easily inequality \eqref{y_0} since it writes
  \begin{multline*}
  -\int_0^LTy_0^2(x)dx+\int_0^T\int_0^Ly^2dxdt +\int_0^T (T-t)y_x^2(0,t)dt \\
-\int_0^T(T-t)(\alpha y_x(0,t)+\beta z(1,t))^2dt=0.
  \end{multline*}
\noindent $\bullet$ Finally, taking $p(\rho,t)=1$ in \eqref{eq_p_z} brings inequality \eqref{z_0} since it writes
$$
\int_0^1 \left(z^2(\rho,T)-z_0^2(-\rho h)\right)d\rho
+\frac{1}{h}\int_0^T  \left(z^2(1,t)-y_x^2(0,t)\right)dt=0.
$$

By density of $\mathcal D(\mathcal A)$ in $H$, the results extend to arbitrary $(y_0,z_0(-h~\cdot))\in H$.
\endproof
 
\subsection{KdV linear equation with a source term 
}
Consider now the KdV linear equation with a right hand side:
\begin{equation}\label{syst_linear_rhs}
\left\{~\begin{array}{l}
y_t(x,t)+y_{xxx}(x,t)+y_x(x,t)=f(x,t), \\
\hfill{~} x\in(0,L),\,t>0,\\
y(0,t)=y(L,t)=0, \hfill{~} t>0,\\
y_x(L,t)=\alpha y_x(0,t)+\beta y_x(0,t-h), \hfill{~} t>0,\\
y_x(0,t)=z_0(t), \hfill{~} t\in(-h,0),\\
y(x,0)=y_0(x), \hfill{~} x\in(0,L).
\end{array}\right.
\end{equation}

\begin{prop}\label{prop:wellposed_rhs}
{\it Assume that \eqref{hyp:alpha_beta} holds.
For any $(y_0,z_0(-h~\cdot))\in H$ and $f\in L^1(0,T,L^2(0,L))$, there exists a unique mild solution $(y,y_x(0,t-h~\cdot))\in \mathcal B \times C([0,T],L^2(0,1))$ to \eqref{syst_linear_rhs}. 
Moreover, denoting $\mathcal E_0 = \left\|(y_0,z_0(-h~\cdot))\right\|^2_{H} + \left\|f\right\|^2_{L^1(0,T,L^2(0,L))} $,  there exists $C>0$ such that 
\begin{eqnarray}
\left\|(y,z)\right\|^2_{C([0,T],H)}\leq C \mathcal E_0 , \label{eq:E_linear_rhs}\\
\| y_x\|^2_{L^2(0,T;L^2(0,L))} \leq C \mathcal E_0 . \label{eq:yH1_rhs}
\end{eqnarray}
}
\end{prop}

\proof
The well-posedness of system \eqref{syst_linear_rhs} in $C([0,T],H)$, when we rewrite it as a first order system (see \eqref{EqA}) with source term $(f(\cdot,t),0)$, stems from $\mathcal A$ being the infinitesimal generator of a $C_0$-semigroup of contractions on $H$. 

Let us now consider a regular solution (which is possible by taking $(y_0,z_0(-h~\cdot))\in\mathcal D(\mathcal A)$ and $f\in L^1([0,T],L^2(0,L))$ for instance) and calculate the time derivative of the energy $E$ defined by \eqref{def:E}:
$$
\begin{aligned}
	\ds \frac{d}{dt}\,E(t) & =  \ds 2\int_0^L y(x,t) y_t(x,t) dx +2\abs{\beta}h\int_0^1 y_x(0,t-h\rho)y_{xt}(0,t-h\rho) d\rho\\
 	& =  \ds -2\int_0^L y(x,t) (y_{xxx}+y_x-f)(x,t) dx \ds - 2 \abs{\beta} \int_0^1 y_x(0,t-h\rho) \partial_\rho y_{x}(0,t-h\rho) d\rho \\
	 & =  \ds y_x^2(L,t)- y_{x}^2(0,t) - \abs{\beta} y_x^2(0,t-h) + \abs{\beta} y_x^2(0,t) + 2\int_0^L f(x,t)y(x,t) dx\\
	 & =  \ds \left(M X(t),X(t)\right) + 2\int_0^L f(x,t)y(x,t) dx,
\end{aligned}
$$
where $X(t)=\begin{bmatrix}y_x(0,t)\\ y_x(0,t-h)\end{bmatrix}$ and $M$ is defined by \eqref{def:M}.
As in the proof of Theorem \ref{thm:wellposed}, under assumption \eqref{hyp:alpha_beta}, $M$ is definite negative, and consequently, there exists $C>0$ such that:
\begin{eqnarray}\label{decay:E}
\dfrac{d}{dt}E(t)\leq -C\left(y_x^2(0,t)+y_x^2(0,t-h)\right)+~2\int_0^L f(x,t)y(x,t) dx\leq 2\|f(t) \|_{L^2(0,L)} \|y(t) \|_{L^2(0,L)}, 
\end{eqnarray}
using also Cauchy-Schwarz inequality.
Integrating between $0$ and $t$, we obtain
$$
\left\|y(.,t)\right\|^2_{L^2(0,L)} + \abs{\beta}h \left\|y_x(0,t-h~\cdot)\right\|^2_{L^2(0,1)}
\leq 
\left\|(y_0,z_0(-h~\cdot))\right\|^2_{H} ~+~ 2\int_0^t \|f(s) \|_{L^2(0,L)} \|y(s) \|_{L^2(0,L)} ds.
$$
Using Young's inequality, for any $\epsilon>0$, we can write
$$	\ds\sup_{t\in[0,T]}\left\|(y(.,t),z(.,t))\right\|^2_{H} \leq  \ds\left\|(y_0,z_0(-h~\cdot))\right\|^2_{H} \\
	+ \epsilon\sup_{t\in[0,T]}\left\|y(.,t)\right\|^2_{L^2(0,L)}+\frac{1}{\epsilon}\left\|f\right\|^2_{L^1(0,T,L^2(0,L))}
$$
and taking $\epsilon$ small enough, then there exists $C>0$ such that
$\ds\sup_{t\in[0,T]}\left\|(y(.,t),z(.t))\right\|^2_{H} \leq C\mathcal E_0, $
yielding \eqref{eq:E_linear_rhs}.

The proof of  \eqref{eq:yH1_rhs} follows exactly the steps of the proof of Proposition \ref{regularity}. One has to pay attention to the right hand side terms that are not homogeneous anymore (but involve the source $f$). Using the first two steps, one proves 
$$
\|y_x(0,.)\|^2_{L^2(0,T)} \leq C\mathcal E_0
$$ and then, \eqref{eq:yH1_rhs} is obtained with the third step.
 
Besides, by density of $\mathcal D(\mathcal A)$ in $H$, the result extends to arbitrary $(y_0,z_0(-h~\cdot))\in H$.
\endproof

\subsection{Well-posedness result of the nonlinear equation}

We endow the space $\mathcal B$ with the norm
$$\left\|y\right\|_{\mathcal B} = \max_{t\in[0,T]} \left\|y(.,t)\right\|_{L^2(0,L)}+\left(\int_0^T\left\|y(.,t)\right\|^2_{H^1(0,L)}dt\right)^{1/2}.
$$

To prove the well-posedness result of the nonlinear system \eqref{syst_nonlinear}, we exactly follow \cite{Coron_2004} (see also \cite{Cerpa_2014}). The proof is given here for the sake of completeness.

The first step is to show that the nonlinearity term $yy_x$ can be considered as a source term of the linear equation \eqref{syst_linear_rhs}:
\begin{prop}\label{prop:Cerpa}
{\it Let $y\in L^2(0,T,H^1(0,L)):=L^2(H^1)$. Then $yy_x\in L^1(0,T,L^2(0,L))$ and the map 
$$y\in L^2(H^1)\mapsto yy_x\in L^1(0,T,L^2(0,L))$$
 is continuous. In particular, there exists $K>0$ such that, for any $y, \tilde y \in L^2(H^1)$, we have
$$\int_0^T \!\!\! \left\|y y_x - \tilde y \tilde y_x \right\|_{L^2(0,L)}
\leq K \left(\left\|y\right\|_{L^2(H^1)}+\left\|\tilde y\right\|_{L^2(H^1)}\right)\left\|y-\tilde y\right\|_{L^2(H^1)}.
$$}
\end{prop}

\proof
The proof can be found in \cite{Rosier_1997} or \cite{Cerpa_2014}.
\endproof\\

Let $(y_0,z_0)\in \mathcal H$ such that $\left\|(y_0,z_0)\right\|_{\mathcal H}\leq r$ where $r>0$ is chosen small enough later. Given $y\in\mathcal B$, we consider the map $\Phi:\mathcal B\rightarrow\mathcal B$ defined by 
$\Phi(y)=\tilde y$ where $\tilde y $ is solution of
$$
\left\{\begin{array}{ll}
\tilde y_t(x,t)+\tilde y_{xxx}(x,t)+\tilde y_x(x,t) = -y(x,t)y_x(x,t), \\
\hfill{} x\in(0,L),\,t>0,\\
\tilde y(0,t)=\tilde y(L,t)=0, \quad t>0,\\
\tilde y_x(L,t)=\alpha \tilde y_x(0,t)+\beta \tilde y_x(0,t-h), \quad t>0,\\
\tilde y_x(0,t)=z_0(t), \quad t\in(-h,0),\\
\tilde y(x,0)=y_0(x), \quad x\in(0,L).
\end{array}\right.
$$
Clearly $y\in\mathcal B$ is a solution of \eqref{syst_nonlinear} if and only if $y$  is a fixed point of the map $\Phi$.
From \eqref{eq:E_linear_rhs}, \eqref{eq:yH1_rhs} and Proposition \ref{prop:Cerpa}, we get
$$\begin{aligned}
&\ds\left\|\Phi(y)\right\|_{\mathcal B} = \left\|\tilde y\right\|_{\mathcal B} 
 \leq  \ds C \left(\left\|(y_0,z_0)\right\|_{\mathcal H}  \! + \!  \int_0^T \! \!  \|y y_x(t)\|_{L^2(0,L)} dt\right)\\
& \leq  \ds C \left(\left\|(y_0,z_0)\right\|_{\mathcal H} + K \|y\|^2_{L^2(0,T,H^1(0,L))}\right) \leq  \ds C \left(\left\|(y_0,z_0)\right\|_{\mathcal H} + \left\|y \right\|^2_{\mathcal B}\right).
\end{aligned}$$
Moreover, for the same reasons, we have
\begin{eqnarray*}
\ds\left\|\Phi(y_1)-\Phi(y_2)\right\|_{\mathcal B} 
& \leq & \ds C \int_0^T\left\|-y_1 y_{1,x}+y_2 y_{2,x}\right\|_{L^2(0,L)}\\
 &\leq&  \ds C  \left(\left\|y_1\right\|_{\mathcal B}+\left\|y_2\right\|_{\mathcal B}\right)\left\|y_1-y_2 \right\|_{\mathcal B}.
\end{eqnarray*}
We consider $\Phi$ restricted to the closed ball $\left\{y\in\mathcal B, \left\|y\right\|_{\mathcal B}\leq R\right\}$ with $R>0$ to be chosen later. 
Then $\left\|\Phi(y)\right\|_{\mathcal B}\leq C \left(r + R^2\right)
$
and
$\left\|\Phi(y_1)-\Phi(y_2)\right\|_{\mathcal B}\leq 2C   R\left\|y_1-y_2 \right\|_{\mathcal B}
$
so that if we take $R$ and $r$ satisfying
$$R<\frac{1}{2C}\qquad\hbox{ and }\qquad r<\frac{R}{2C},
$$
then
$\left\|\Phi(y)\right\|_{\mathcal B}<R$ and $\|\Phi(y_1)-\Phi(y_2)\|_{\mathcal B}\leq 2 C R\|y_1-y_2 \|_{\mathcal B}$,
with $2C R<1$. Consequently, we can apply the Banach fixed point theorem and the map $\Phi$ has a unique fixed point. We have then shown the following:
\begin{prop}
{\it Let $T,L>0$ and assume that \eqref{hyp:alpha_beta} holds. Then there exist $r>0$ and $C>0$ such that for every $(y_0,z_0)\in \mathcal H$ verifying 
$\left\|(y_0,z_0)\right\|_{\mathcal H}\leq r,$
there exists a unique   $y \in \mathcal B$ solution of system \eqref{syst_nonlinear} which satisfies
$
\left\|y\right\|_{\mathcal B}\leq C \left\|(y_0,z_0)\right\|_{\mathcal H}.
$
}
\end{prop}

\section{Lyapunov approach for a first stabilization result}\label{LyapSec}
The goal of this section is to prove our first main result, presented in Theorem~\ref{thm:expostab_nonlinear}.  We will basically detail the proof of the exponential stability of the solution of system~\eqref{syst_nonlinear}, which is based on the appropriate choice of a candidate Lyapunov functional. A first step is the following proposition concerning the energy of the system.

\begin{prop}\label{prop:E'_nl}
Let \eqref{hyp:alpha_beta} be satisfied.
Then, for any regular solution of \eqref{syst_nonlinear} the energy $E$ defined by \eqref{def:E} is non-increasing and satisfies
\begin{equation}\label{decay:E_nl}
\dfrac{d}{dt}\,E(t)= \left(\alpha^2-1+\abs{\beta}\right) y_x^2(0,t)+\left(\beta^2-\abs{\beta}\right)y_x^2(0,t-h) 
+2\alpha\beta y_x(0,t)y_x(0,t-h)\leq 0.
\end{equation}
\end{prop}

\proof
Differentiating \eqref{def:E} and using \eqref{syst_nonlinear}, we obtain 
$$\begin{aligned}
\ds \dfrac{d}{dt}\,E(t)  =&  
   \ds -2\int_0^L y(x,t)( y_{xxx}+y_x+yy_x)(x,t) dx - 2 \abs{\beta} \int_0^1 y_x(0,t-h\rho) \partial_\rho y_{x}(0,t-h\rho) d\rho \\
 =&~  \ds y_x^2(L,t)- y_{x}^2(0,t) - \abs{\beta} y_x^2(0,t-h) + \abs{\beta} y_x^2(0,t)\\
 =& ~ \left(\alpha^2-1+\abs{\beta}\right) y_x^2(0,t)+\left(\beta^2-\abs{\beta}\right)y_x^2(0,t-h) +2\alpha\beta y_x(0,t)y_x(0,t-h)\\
 =&~  \ds \left( M X(t),X(t)\right),
\end{aligned}$$
where as usual $X(t)=\left(y_x(0,t), y_x(0,t-h)\right)$ and $M$ is defined by \eqref{def:M}.
As in the proof of Theorem~\ref{thm:wellposed}, under assumption \eqref{hyp:alpha_beta}, $M$ is definite negative, ending  the proof.
\endproof

\begin{rk}
We deduce from assumptions \eqref{hyp:alpha_beta}  and Proposition~\ref{prop:E'_nl} that the energy $E$ is decreasing as long as $y_x(0, t)$ does not vanish.\\
\end{rk}
This result on the energy of the system does not yield the exponential stability we are seeking. Therefore, we choose now the following candidate Lyapunov functionnal:
\begin{equation}\label{def:Lyapunov}
V(t) = E(t) + \mu_1 V_1(t) + \mu_2 V_2(t), 
\end{equation}
where $\mu_1$ and $\mu_2\in(0,1)$ are positive constants that will be fixed small enough later on, $E$ is the energy defined by \eqref{def:E}, $V_1$ is defined by 
\begin{equation}\label{def:V1}
V_1(t) = \int_0^L x y^2(x,t) dx,
\end{equation}
and
$V_2$ is defined by 
\begin{equation}\label{def:V2}
V_2(t) = h \int_0^1 (1-\rho) y_x^2(0,t-h\rho) d\rho.
\end{equation}
It is clear that the two energies $E$ and $V$ are equivalent, in the sense that
\begin{equation}\label{eq:EV_equivalent}
E(t) \leq V(t) \leq \left(1+\max\left\{L\mu_1, \frac{\mu_2}{\abs{\beta}}\right\}\right) E(t).
\end{equation}

\medskip\goodbreak\noindent {\textbf{Proof of Theorem~\ref{thm:expostab_nonlinear}.}\\
Let $y$ be a regular solution of \eqref{syst_nonlinear}.
Differentiating \eqref{def:V1} and using \eqref{syst_nonlinear}, we obtain 
$$\begin{aligned}
&\ds \frac{d}{dt}\, V_1(t)   =  \ds -2\int_0^L xy(x,t) (y_{xxx}+y_x+yy_x)(x,t) dx \\
 & =  \ds -2\int_0^L y_x^2(x,t) dx + 2\left[y(x,t) y_{x}(x,t)\right]_0^L - \int_0^L y_x^2(x,t)dx \\
 &+\left[x y_x^2(x,t)\right]_0^L   \ds + \int_0^L y^2(x,t) dx + \frac{2}{3} \int_0^L y^3(x,t)dx\\
 & =  \ds -3 \int_0^L y_x^2(x,t)dx + L\alpha^2 y_x^2(0,t)+L\beta^2 y_x^2(0,t-h) +2\alpha \beta L y_x(0,t) y_x(0,t-h)\\
& \ds  + \int_0^L y^2(x,t) dx + \frac{2}{3} \int_0^L y^3(x,t)dx.
\end{aligned}$$
Moreover, differentiating \eqref{def:V2} and using \eqref{syst_nonlinear}, we obtain 
$$\begin{aligned}
&\ds \frac{d}{dt}\, V_2(t)   =  \ds -2 \int_0^1 (1-\rho) y_x(0,t-h\rho) \partial_\rho y_{x}(0,t-h\rho) d\rho 
  =  \ds y_x^2(0,t) - \int_0^1  y_{x}^2(0,t-h\rho) d\rho.
\end{aligned}$$
Consequently, for any $\gamma>0$, we have 
$$\begin{aligned}
&\ds \frac{d}{dt}\, V(t) + 2\gamma V(t)  =  \ds\left(\alpha^2-1+\abs{\beta}+\mu_1 L\alpha^2+\mu_2\right) y_x^2(0,t)\\
&+\left(\beta^2-\abs{\beta}+\mu_1 L\beta^2\right)y_x^2(0,t-h) 
  \ds + 2\alpha\beta\left(1+L\mu_1\right)y_x(0,t)y_x(0,t-h)-3\mu_1 \int_0^L \!\! y_x^2(x,t)dx \\
& \ds +\left(2\gamma \abs{\beta}h-\mu_2\right)  \!\! \int_0^1 \!  \!\! y_{x}^2(0,t-h\rho) d\rho +(2\gamma+\mu_1)\int_0^L \! \!\! y^2(x,t) dx  \\
&  \ds + 2\gamma \mu_1 \int_0^L x y^2(x,t) dx + 2\gamma\mu_2 h \int_0^1 (1-\rho) y_x^2(0,t-h\rho) d\rho  + \frac{2}{3}\,\mu_1 \int_0^L y^3(x,t)dx\\
& \leq \ds \left(M_{\mu_1}^{\mu_2} X(t),X(t)\right)-3\mu_1 \int_0^L y_x^2(x,t)dx 
+\left(2\gamma h(\mu_2+\abs{\beta})-\mu_2\right)\int_0^1  y_{x}^2(0,t-h\rho) d\rho\\
&  \ds + \left(2\gamma\left(1+L\mu_1\right)+\mu_1\right)\int_0^L y^2(x,t) dx + \frac{2}{3}\,\mu_1 \int_0^L y^3(x,t)dx,
\end{aligned}$$
where 
$
X(t)=\begin{bmatrix}y_x(0,t)\\ y_x(0,t-h)\end{bmatrix}
$ 
 and 
$
M_{\mu_1}^{\mu_2}=\begin{bmatrix}(1+L\mu_1)\alpha^2-1+\abs{\beta}+\mu_2 & \alpha\beta\left(1+L\mu_1\right) \\ 
\alpha\beta\left(1+L\mu_1\right) & (1+L\mu_1)\beta^2-\abs{\beta}\end{bmatrix},
$ noticing  that 
$$
M_{\mu_1}^{\mu_2}= M+\mu_1 L\begin{pmatrix}\alpha^2&\alpha\beta\\\alpha\beta & \beta^2\end{pmatrix}+\mu_2\begin{pmatrix}1&0\\0&0\end{pmatrix},
$$
where $M$ is defined by \eqref{def:M}.\\
As $M$ is definite negative, we easily prove that for $\mu_1$ and $\mu_2>0$ sufficiently small the matrix $M_{\mu_1}^{\mu_2}$ is definite negative, by continuity of the applications Trace and Determinant.

Finally, for $\mu_1$ and $\mu_2$ sufficiently small, using Poincar\'e inequality, ($ \|y \|_{L^2(0,L)}\leq \frac{L}{\pi} \|y_x \|_{L^2(0,L)}$ for $y\in H^1_0(0,L)$) we obtain 
\begin{multline*}
\frac{d}{dt}\, V(t)+ 2\gamma V(t)
\leq 
\left(\frac{L^2\left(2\gamma\left(1+L\mu_1\right)+\mu_1\right)}{\pi^2}-3\mu_1\right)  \|y_x(t)\|^2_{L^2(0,L)}  \\
+\left(2\gamma h(\mu_2+\abs{\beta})-\mu_2\right)  \int_0^1  y_{x}^2(0,t-h\rho) d\rho + \frac{2}{3}\,\mu_1 \int_0^L y^3(x,t)dx.
\end{multline*}
Moreover, using Cauchy-Schwarz inequality, Proposition~\ref{prop:E'_nl} and since $H^1_0(0,L)\ \subset L^{\infty}(0,L)$, we have:
$$\begin{array}{rcl}
\ds \int_0^L y^3(x,t)dx & \leq & \ds \left\|y(.,t)\right\|^2_{L^{\infty}(0,L)}\int_0^L |y(x,t)|dx\\
 & \leq & \ds L \sqrt{L} \left\|y_x(.,t)\right\|^2_{L^2(0,L)} \left\|y(.,t)\right\|_{L^{2}(0,L)}\\
 & \leq & \ds L^{3/2} \left\|(y_0,z_0)\right\|_{\mathcal H}\left\|y_x(.,t)\right\|^2_{L^{2}(0,L)}\\
 & \leq & \ds L^{3/2}\, r \left\|y_x(.,t)\right\|^2_{L^{2}(0,L)}.
\end{array}$$
Consequently, we have
$$\frac{d}{dt}\, V(t)+ 2\gamma V(t) \leq \Upsilon  \|y_x(t)\|^2_{L^2(0,L)} 
+\left(2\gamma h(\mu_2+\abs{\beta})-\mu_2\right)\int_0^1  y_{x}^2(0,t-h\rho) d\rho
$$
where $\Upsilon = \dfrac{L^2\left(2\gamma\left(1+L\mu_1\right)+\mu_1\right)}{\pi^2}-3\mu_1+\dfrac{2L^{3/2}r\mu_1}{3}.$\\
Since $L$ satisfies the constraint \eqref{hyp:L}, it  is possible to  choose $r$  small enough to have
$r<\dfrac{3(3\pi^2-L^2)}{2L^{3/2}\pi^2}$. Then one can choose $\gamma>0$ such that \eqref{eq:gamma} holds in order to obtain
$$\frac{d}{dt}\, V(t)+ 2\gamma V(t)\leq 0, ~~\forall t >0.$$
Integrating over $(0,t)$ and using \eqref{eq:EV_equivalent}, we finally obtain that
$$E(t)\leq \left(1+\max\left\{L\mu_1, \frac{\mu_2}{\abs{\beta}}\right\}\right) E(0)e^{-2\gamma t},\qquad \forall t>0.
$$
By density of $\mathcal D(\mathcal A)$ in $H$, the results extend to arbitrary $(y_0,z_0)\in\mathcal H$.
\endproof\\

\begin{rk} {\bf On the size of the delay.}
As one can deduce from \eqref{eq:gamma}, when the delay $h$ increases, the decay rate $\gamma$ decreases.\\
\end{rk}

\begin{rk} {\bf On the coefficients of the boundary feedback.} In the case without delay (i.e. $\beta=0$), adapting the previous results, the exponential stability of Theorem \ref{thm:expostab_nonlinear} is satisfied if and only if
	$\left|\alpha\right|<1$,
	which corresponds to the assumptions given in \cite{Zhang_1994} and \cite{Perla_2002}.
	 In the case where $\alpha=0$, the exponential stability of Theorem \ref{thm:expostab_nonlinear} is satisfied if and only if $\left|\beta\right|<1$.
	 Even if $\alpha=\beta=0$, Theorem \ref{thm:expostab_nonlinear} shows that the system \eqref{syst_nonlinear} is exponentially stable in the case where $0<L<\sqrt 3\pi$ (see \cite{Perla_2002} in the case where $L\notin\mathcal N$). We recall that the main goal of this paper is to show that a delay does not destabilize the system, which may be the case in many other delayed systems (see, for instance, \cite{Datko_1988}, \cite{Datko_1986}, \cite{Logemann_1996}, \cite{Rebarber_1999}). \\
\end{rk}

\begin{rk}\label{rk:mubounded}
{\bf On the  estimation of  parameters $r$, $\mu_1$, $\mu_2$.} By the proof of Theorem \ref{thm:expostab_nonlinear}, $\mu_1$ and $\mu_2$ are chosen small enough such that $M_{\mu_1}^{\mu_2}$ is definite negative. Simple calculations show that taking $\mu_1>0$ and $\mu_2\in(0,1)$ such that
$$\mu_2<\min\left\{1-\alpha^2-\beta^2, \frac{(\left|\beta\right|-1)^2-\alpha^2}{1-\left|\beta\right|},\frac{\alpha^2-\beta^2+\left|\beta\right|}{\left|\beta\right|}\right\},$$
$$
\mu_1<\min\left\{\frac{1-\mu_2-(\alpha^2+\beta^2)}{L(\alpha^2+\beta^2)}, \right.
\left.\frac{(\left|\beta\right|-1)^2-\alpha^2-\mu_2(1-\left|\beta\right|)}{L(\alpha^2-\beta^2+\left|\beta\right|(1-\mu_2))}\right\}
$$
implies that $M_{\mu_1}^{\mu_2}$ is definite negative. It is then sufficient to take $r>0$ such that
$r<\frac{3(3\pi^2-L^2)}{2L^{3/2}\pi^2}$
to have the exponential decreasing of the energy $E$ with a decay rate $\gamma$ given by \eqref{eq:gamma}. \\
\end{rk}

\begin{rk} {\bf On the stabilization of the linear KdV equation.}
With the same Lyapunov functional \eqref{def:Lyapunov}, we can obtain the exponential stability result for the linear equation \eqref{syst_linear} under the assumptions \eqref{hyp:alpha_beta} and \eqref{hyp:L} (and without restriction about the initial data) for $\mu_1$ and $\mu_2$ small enough. In this case the decay rate is given by
$$
\gamma\leq\min\left\{\frac{(3\pi^2-L^2)\mu_1}{2L^2(1+L\mu_1)},\frac{\mu_2}{2(\mu_2+\abs{\beta})h}\right\}. 
$$ 
\end{rk}

\begin{rk} {\bf On the length of the spacial domain.}
Our hypothesis $L<\sqrt 3\pi$ (see \eqref{hyp:L}) eliminates the set of critical lengths $\mathcal N$ (see Rosier \cite{Rosier_1997}), but also lots of non critical lengths. The condition $L<\sqrt 3\pi$ is a technical one and comes from the choice of the multiplier $x$ in the expression of $V_1$. To find a better multiplier is an open problem as far as we know. \\
\end{rk}

In the next section we will prove a stabilization result for all non critical lengths but without any bound on the decay rate.

\section{Second stabilization result - Observability approach}\label{ObsSec}

This section aims at proving our second main result, stated in Theorem~\ref{thm:expostab_nonlinear_LN}, which is obtained simply for non critical lengths and gives generic exponential stability of the solution of system~\eqref{syst_nonlinear}.
The proof relies on an observability inequality and the use of a contradiction argument.  It will need several steps in order to handle the nonlinearity of the KdV equation under consideration.

 \subsection{Proof of the stability of the linear equation}

We first prove the following observability result.
\begin{theo}\label{observability}
{\it Assume that \eqref{hyp:alpha_beta} is satisfied.
Let $L\in (0,+\infty)\setminus \mathcal N$ and $T>h$. Then there exists $C>0$ such that for all $(y_0,z_0(-h~\cdot))\in H$, 
we have the observability inequality
\begin{equation}\label{ineq_observ}
\int_0^L y^2_0(x)dx+\abs{\beta}h\int_0^1z_0^2(-h\rho)d\rho
\leq C\int_0^T\left(y_x^2(0,t)+z^2(1,t)\right)dt
\end{equation}
where $(y,z)=S(.)(y_0,z_0(-h~\cdot))$.
}
\end{theo}

\proof We proceed by contradiction as in \cite{Rosier_1997} . Let us suppose that \eqref{ineq_observ} is false. Then there exists a sequence $\bigg((y_0^n,z_0^n(-h~\cdot))\bigg)_n\subset H$ such that 
$$\int_0^L (y^n_0)^2(x)dx+\abs{\beta}h\int_0^1(z^n_0)^2(-h\rho) d\rho = 1$$ 
and 
\begin{equation}\label{eq:contradiction2}
\norm{y^n_x(0,.)}^2_{L^2(0,T)}+\norm{z^n(1,.)}^2_{L^2(0,T)}\rightarrow 0\hbox{ as }n\rightarrow +\infty,
\end{equation}    
where $(y^n,z^n)=S(y_0^n,z_0^n(-h~\cdot))$. Thanks to Proposition~\ref{regularity}, $(y^n)_n$ is a bounded sequence in $L^2(0,T,H^1(0,L))$, and then $y^n_t=-y^n_x-y^n_{xxx}$ is bounded in $L^2(0,T,H^{-2}(0,L))$. Due to a result of Simon \cite{Sim87}, the set $\{y^n\}_n$ is relatively compact in $L^2(0,T,L^2(0,L))$ and we may assume that $(y^n)_n$ is convergent in $L^2(0,T,L^2(0,L))$. 

Thanks to \eqref{y_0} and \eqref{eq:contradiction2}, we deduce that $(y_0^n)_n$ is a Cauchy sequence in $L^2(0,L)$. 
We now prove that if $T>h$, $(z_0^n(-h~\cdot))_n$ is a Cauchy sequence in $L^2(0,1)$. Indeed, as $z^n(\rho,T)=y_x^n(0,T-\rho h)$, if $T>h$, we have
$$
\ds \int_0^1 (z^n(\rho,T))^2d\rho  =  \ds \int_0^1 (y_x^n(0,T-\rho h))^2d\rho \leq  \ds \frac{1}{h}\int_0^T (y_x^n(0,t))^2dt.
$$
Using \eqref{z_0}, for $T>h$ we have
$$
\norm{z_0^n(-h~\cdot)}^2_{L^2(0,1)}\leq \frac{1}{h}\norm{y_x^n(0,.)}^2_{L^2(0,T)}+\frac{1}{h}\norm{z^n(1,.)}^2_{L^2(0,T)}.
$$
Thus $(z_0^n(-h~\cdot))_n$ is a Cauchy sequence in $L^2(0,1)$ using also \eqref{eq:contradiction2}.

Let $(y_0,z_0(-h~\cdot))=\lim (y_0^n,z_0^n(-h~\cdot))$ in $H$ and $(y,z)=S(.)(y_0,z_0(-h~\cdot))$. By using Proposition~\ref{regularity}, $(y_x(0,.), z(1,.))=\lim (y_x^n(0,.),z^n(1,.))$ in $L^2(0,T)$. Thus we have that $\int_0^L y^2_0(x)dx+\abs{\beta}h\int_0^1z_0^2(-h\rho)d\rho = 1$ and $(y_x(0,.), z(1,.))=0$.
As $z(1,t)=y_x(0,t-h)=0$ we deduce that $z_0=0$ and $z=0$. 
Consequently $y$ is solution of
$$
\left\{\begin{array}{ll}
y_t(x,t)+y_{xxx}(x,t)+y_x(x,t)=0, \hfill{~} x\in(0,L),\,t>0,\\
y(0,t)=y(L,t)=0, \hfill{~} t>0,\\
y_x(L,t)=y_x(0,t)=0, \hfill{~} t>0,\\
y(x,0)=y_0(x), \hfill{~} x\in(0,L),
\end{array}\right.
$$
and $\left\|y_0\right\|_{L^2(0,L)}=1$.\\
We can apply the result of \cite[Lemma 3.4]{Rosier_1997}: if $L\notin \cal N$ there is no function satisfying this last system. Then we obtain a contradiction, which ends the proof of 
Theorem \ref{observability}.
\endproof\\

From observability inequality \eqref{ineq_observ}, one can deduce the exponential stability of the KdV linear system \eqref{syst_linear}, stated here:
\begin{theo}\label{thm:expostab_linear_N}
{\it Assume that $L\in (0,+\infty)\setminus \mathcal N$ and that \eqref{hyp:alpha_beta} is satisfied.
Then, for every $(y_0,z_0)\in\mathcal H$, the energy of system \eqref{syst_linear}, denoted by $E$ and defined by \eqref{def:E}, decays exponentially. More precisely, there exist two positive constants $\nu$ and $\kappa$ such that 
$E(t)\leq \kappa E(0)e^{-\nu t}$, for all $t>0$.
} 
\end{theo}

\proof
We follow the same kind of proof as in \cite{Nicaise_2006}. Let $(y_0,z_0)\in\mathcal D(\mathcal A)$. Integrating \eqref{decay:E} with $f=0$ between $0$ and $T>h$, we have
$$E(T)-E(0)\leq -C_1 \int_0^T \left(y_x^2(0,t)+y_x^2(0,t-h)\right) dt,
$$
which is equivalent to
\begin{equation}\label{eq:E'integr}
\hspace{-0.3cm}\int_0^T \left(y_x^2(0,t)+y_x^2(0,t-h)\right) dt \leq \frac{1}{C_1}\left(E(0)-E(T)\right).
\end{equation}
As the energy is non-increasing, we have, using the observability inequality \eqref{ineq_observ} and \eqref{eq:E'integr},
$$
E(T)\leq E(0) \leq C\int_0^T \left(y_x^2(0,t)+y_x^2(0,t-h)\right) dt 
\leq\frac{C}{C_1}\left(E(0)-E(T)\right),
$$
which implies that
\begin{equation}\label{eq:decayE_gamma}
E(T)\leq \gamma E(0), \hbox{ with } \gamma=\frac{\frac{C}{C_1}}{1+\frac{C}{C_1}}<1.
\end{equation}
Using  this argument on $[(m-1)T,mT]$ for $m=1,2,...$ (which is valid because the system is invariant by translation in time), we will get 
$$E(mT)\leq \gamma E((m-1)T)\leq\cdots\leq \gamma^m E(0).
$$
Therefore, we have
$E(mT)\leq e^{-\nu mT }E(0)$
with $\nu=\frac{1}{T}\ln\frac{1}{\gamma}=\frac{1}{T}\ln\left(1+\frac{C_1}{C}\right)>0$. For an arbitrary positive $t$, there exists $m\in\mathbb N^*$ such that $(m-1)T<t\leq mT$, and by the non-increasing property of the energy, we conclude that
$$E(t)\leq E((m-1)T) \leq e^{-\nu (m-1)T }E(0)\leq\frac{1}{\gamma} e^{-\nu t}E(0).
$$ 
By density of $\mathcal D(\mathcal A)$ in $H$, we deduce that the exponential decay of the energy $E$ holds for any initial data in $\mathcal H$.
\endproof

\subsection{Stability of the nonlinear equation}

We consider in this section the more general case than in Theorem \ref{thm:expostab_nonlinear}, where $L\in(0,+\infty)\backslash\mathcal N$, and prove the exponential decay of small amplitude solutions of the nonlinear KdV equation \eqref{syst_nonlinear}:

\medskip\goodbreak\noindent {\textbf{Proof of Theorem~\ref{thm:expostab_nonlinear_LN}.}\\
The proof follows \cite{Cerpa_2014} for the stabilization of the nonlinear KdV equation with internal feedback.
Consider initial data $\left\|(y_0,z_0)\right\|_{\mathcal H}\leq r$ with $r$ chosen later. The solution $y$ of \eqref{syst_nonlinear} can be written as $y=y^1+y^2$ where $y^1$ is solution of
$$\left\{\begin{array}{ll}
y^1_t(x,t)+y^1_{xxx}(x,t)+y^1_x(x,t)=0, & x\in(0,L),\,t>0,\\
y^1(0,t)=y^1(L,t)=0, & t>0,\\
y^1_x(L,t)=\alpha y^1_x(0,t)+\beta y^1_x(0,t-h), & t>0,\\
y^1_x(0,t)=z_0(t), & t\in(-h,0),\\
y^1(x,0)=y_0(x), & x\in(0,L),
\end{array}\right.
$$
and $y^2$ is solution of
$$\left\{\begin{array}{ll}
y^2_t(x,t)+y^2_{xxx}(x,t)+y^2_x(x,t)=-y(x,t)y_x(x,t), \\
\hfill{} x\in(0,L),\,t>0,\\
y^2(0,t)=y^2(L,t)=0, \qquad t>0,\\
y^2_x(L,t)=\alpha y^2_x(0,t)+\beta y^2_x(0,t-h), \qquad t>0,\\
y^2_x(0,t)=0, \qquad t\in(-h,0),\\
y^2(x,0)=0, \qquad x\in(0,L).
\end{array}\right.
$$
More precisely, $y^1$ is solution of \eqref{syst_linear} with initial data $(y_0,z_0(-h~\cdot))\in H$ and $y^2$ is solution of \eqref{syst_linear_rhs} with initial data $(0,0)$ and right-hand side $f=-yy_x\in L^1(0,T,L^2(0,L))$ (see Proposition\ref{prop:Cerpa}).
Using \eqref{eq:decayE_gamma}, Proposition~\ref{prop:wellposed_rhs} and Proposition~\ref{prop:Cerpa}, we have
\begin{equation}\label{eq:maj1_nl}
\begin{aligned}
&\left\|(y(T),z(T))\right\|_H  \leq  \left\|(y^1(T),z^1(T))\right\|_H + \left\|(y^2(T),z^2(T))\right\|_H\\
& \leq  \gamma \left\|(y_0,z_0(-h~\cdot))\right\|_H + C \left\|y y_x\right\|_{L^1(0,T,L^2(0,L))}\\
& \leq  \gamma \left\|(y_0,z_0(-h~\cdot))\right\|_H + C \left\|y \right\|^2_{L^2(0,T,H^1(0,L))},
\end{aligned}\end{equation}
with $0<\gamma<1$.
The aim is now to deal with the last term of the previous inequality.
For that, we multiply the first equation of \eqref{syst_nonlinear} by $xy$ and integrate to obtain
\begin{multline*}
3\int_0^T \int_0^L y_x^2(x,t) dx dt + \int_0^L x y^2(x,T) dx \\
= \int_0^T \int_0^L y^2(x,t) dx dt + \int_0^L x y_0^2(x) dx 
+ L \int_0^T y_x^2(L,t) dt + \frac{2}{3} \int_0^T \int_0^L y^3(x,t) dx dt.
\end{multline*}
Consequently, by \eqref{syst_nonlinear} and \eqref{decay:E_nl}, we have
\begin{multline*}
\int_0^T \int_0^L y_x^2(x,t) dx dt \leq \frac{T+L}{3}\left\|(y_0,z_0)\right\|_{\mathcal H}^2 \\
+ \frac{L}{3} \int_0^T \left(\alpha y_x(0,t)+\beta z(1,t)\right)^2 dt + \frac{2}{9} \int_0^T \int_0^L |y|^3(x,t) dx dt.
\end{multline*}
As $H^1(0,L)$ embeds into $C([0,L])$ and using Cauchy-Schwarz inequality and \eqref{decay:E_nl}, we have
$$\begin{aligned}
&\ds \int_0^T \int_0^L |y|^3(x,t) dx dt  \leq  \ds \int_0^T \left\|y\right\|_{L^\infty(0,L)}\int_0^L y^2(x,t) dx dt\\
& \leq  \ds \sqrt{L}\int_0^T \left\|y\right\|_{H^1(0,L)}\int_0^L y^2(x,t) dx dt\\
& \leq  \ds \sqrt{L}\left\|y\right\|^2_{L^\infty(0,T,L^2(0,L))}\int_0^T \left\|y\right\|_{H^1(0,L)} dt
 \leq  \ds \sqrt{LT} \left\|(y_0,z_0)\right\|^2_{\mathcal H} \left\|y\right\|_{L^2(0,T,H^1(0,L))}.
\end{aligned}$$
We deduce that
\begin{multline}\label{eq:a}
\int_0^T \int_0^L y_x^2(x,t) dx dt \leq \frac{T+L}{3}\left\|(y_0,z_0)\right\|_{\mathcal H}^2 
+ \frac{L}{3} \int_0^T \left(\alpha y_x(0,t)+\beta z(1,t)\right)^2 dt \\
+ \frac{2\sqrt{LT}}{9} \left\|(y_0,z_0)\right\|^2_{\mathcal H} \left\|y\right\|_{L^2(0,T,H^1(0,L))}.
\end{multline}
Now we multiply the first equation of \eqref{syst_nonlinear} by $y$ and integrate to obtain
$$\int_0^L y^2(x,T)dx - \int_0^L y_0^2(x)dx\\
-\int_0^T(\alpha y_x(0,t)+\beta z(1,t))^2dt
+\int_0^T y_x^2(0,t)dt=0.
$$
Using the same proof as in the linear case (see \eqref{ineq_z} e.g.), we obtain that 
$$\int_0^T y_x^2(0,t) dt + \int_0^T z^2(1,t) dt \leq C\left\|(y_0,z_0)\right\|^2_{\mathcal H}.
$$
Consequently, we have
$$
\int_0^T(\alpha y_x(0,t)+\beta z(1,t))^2dt \leq 2(\alpha^2+\beta^2)\int_0^T(y_x^2(0,t)
+z^2(1,t))dt
\leq 2C(\alpha^2+\beta^2)\left\|(y_0,z_0)\right\|^2_{\mathcal H},
$$
and then, using Young's inequality in \eqref{eq:a}, there exists $C>0$ such that 
\begin{equation}\label{eq:maj2_nl}
\hspace{-0.3cm}\int_0^T\! \! \int_0^L y_x^2(x,t) dx dt \leq C\left(\left\|(y_0,z_0)\right\|_{\mathcal H}^2 + \left\|(y_0,z_0)\right\|^4_{\mathcal H}\right).
\end{equation}
Therefore, gathering \eqref{eq:maj1_nl} and \eqref{eq:maj2_nl}, there exists $C>0$ such that
\begin{equation}
\left\|(y(T),z(T))\right\|_{H}
 \leq  \left\|(y_0,z_0)\right\|_{\mathcal H} \! \left(\gamma + C \left\|(y_0,z_0)\right\|_{\mathcal H} + C \left\|(y_0,z_0)\right\|_{\mathcal H}^3\right)
\end{equation}
which implies 
$$
\left\|(y(T),z(T))\right\|_H \leq  \left\|(y_0,z_0)\right\|_{\mathcal H}\left(\gamma + C r + C r^3\right).
$$
Given $\epsilon>0$ small enough such that $\gamma+\epsilon<1$, we can take $r$ small enough such that $r+r^3<\frac{\epsilon}{C}$, in order to have
$$\left\|(y(T),z(T))\right\|_H \leq  (\gamma+\epsilon)\left\|(y_0,z_0)\right\|_{\mathcal H},$$
with $\gamma+\epsilon<1$. Then end of the proof follows the same lines as in the linear case (see after \eqref{eq:decayE_gamma}).
\endproof

\section{Numerical simulations}\label{NumSec}

In this last section, we propose to illustrate the obtained results with some numerical simulations.  For doing this, we have used the numerical scheme proposed by Colin and Gisclon \cite{colin2001initial}. We first choose a time-step $\delta t$ and a space-step $\delta x$ and we denote by $y_j^n$ the approximate value of $y(j\delta x, n\delta t)$. We also choose a delay-step, $\delta \rho$ and denote by  $z_i^n$ the approximate value of $z(i\delta \rho, n\delta t)$.
Then, we define $J=L/\delta x$ the number of space steps , $K=h/\delta \rho$ the number of  delay steps.

 We denote by $X_D$ the following space of finite sequence,
$$
X_D=\Big\{(y,z)=(y_0,\dots, y_J,z_0,\dots,z_K)\in \RR^{J+K+2}, 
\hbox{ with } y_0=y_J=0\Big\}.
$$
We introduce the classical difference operators
$$(D^+_{\delta x}y)_j=\frac{y_{j+1}-y_j}{\delta x} ,\quad  (D^-_{\delta x}y)_j=\frac{y_{j}-y_{j-1}}{\delta x}$$
and  $ D_{\delta x}=\dfrac{D^+_{\delta x} +D^-_{\delta x}}{2}.$
Then we propose to approach the linear system \eqref{EqA}, by the following implicit numerical scheme,
$$
\frac{(y,z)^{n+1}-(y,z)^n}{\delta t}+A(y,z)^{n+1}=0
$$
where the matrix $A\in \RR^{J-2+K}\times \RR^{J-2+K}$ is an approximation of the operator 
$$-\mathcal A=\begin{pmatrix}\partial_{xxx}+\partial_x& 0\\ 0&\frac{1}{h}\,\partial_{\rho}\end{pmatrix},$$ choosing $\frac{1}{h}D^-_{\delta \rho}$ and  $D^+_{\delta x}D^+_{\delta x}D^-_{\delta x}+D_{\delta x}$ to approximate  $\frac{1}{h}\,\partial_{\rho}$ and $\partial_{xxx}+\partial_x $.\\

In order to understand why $A$ is of size $J-2+K$ instead of $J+K+2$, one can notice that besides $y(0,t)=y(L,t)=0$ implying $y_0^n=0$ and $y_J^n=0$,
we also have $z(0,t)=y_x(0,t)$ so that $z_0^n=\dfrac{y_1^n}{\delta x}$. Furthermore, as $y_x(L,t)=\alpha y_x(0,t)+\beta z(1,t)$ we get as an approximation,
$$ \dfrac{y^n_{J}-y^n_{J-1}}{\delta x}=\alpha \dfrac{y^n_{1}-y^n_0}{\delta x}+\beta z^n_K $$ yielding the following equation:
\begin{equation}
y^n_{J-1}=-\alpha y^n_{1}-\beta \delta x z^n_K.
\end{equation}
This equation, at order $n$, gives the last un-necessary  coefficient, and written at orden $n+1$, 
using also the implicit scheme's equations, gives the terms in $\alpha$ and $\beta$ appearing in the coming matrices $\Lambda_{11}$ and $\Lambda_{12}$.

Let us finally detail the shape of the matrix  ${\bf A}$  defined such that the numerical problem to solve is to find $\{(y,z)^n\}\in X_D$ verifying 
$$
\left\{\begin{array}{l}
{\bf A}(y,z)^{n+1}_{(1:J-2),(1:K)}=(y,z)^{n}_{(1:J-2),(1:K)},\\
z_0^n=\frac{y_1^n}{\delta x}.
\end{array}\right.$$
It writes
${\bf A} = \left(\begin{array}{c|c} \Lambda_{11} & \Lambda_{12}  \\\hline \Lambda_{21} & \Lambda_{22} \end{array}\right)$
where ${\bf 0}$ denotes vectors or matrices of $0$ and of appropriate dimensions and 
\begin{equation*}
\Lambda_{11}=\left (\begin{array}{cccccc}
	a_1&a_2 & a_3 & 0 &\dots &0\\
 	a_4 &  a_1 &   &   & \ddots  & \vdots \\
	0&  & \ddots&  &  & 0  \\
	{\bf 0}& \ddots & &  &  &  a_3  \\
	-\dfrac{\alpha\delta t}{\delta x^3}&{\bf 0}& 0 &a_4  & a_1 & a_2 \\
	\dfrac{3\alpha\delta t}{\delta x^3}-\dfrac{\alpha\delta t}{2\delta x} & 0 & \dots & 0 &a_4 &  a_1   \end{array}
   \right)
\end{equation*}
with $a_1 = 1+\dfrac{3\delta t}{\delta x^3}$, $a_2 = \dfrac{-3\delta t}{\delta x^3}+\dfrac{\delta t}{2\delta x}$, 
$a_3 = \dfrac{\delta t}{\delta x^3}$ and  $a_4 =  \dfrac{-\delta t}{\delta x^3}-\dfrac{\delta t}{2\delta x}$,\\
$
\Lambda_{12}=\left (\begin{array}{cc}
	{\bf 0}&{\bf 0}\\
	{\bf 0}&-\dfrac{\beta\delta t}{\delta x^2}\\
{\bf 0}&\dfrac{3\beta\delta t}{\delta x^2}-\dfrac{\beta\delta t}{2}\\
   \end{array}
   \right)
$, ~
$
\Lambda_{21}=\left (\begin{array}{cccc}
	\dfrac{-\delta t}{h\delta x\delta \rho}&{\bf 0}\\
	{\bf 0}&{\bf 0}
   \end{array}
   \right)
$
and 
$
\Lambda_{22}=\left (\begin{array}{ccc}
	1+\dfrac{\delta t}{h\delta \rho}&0&{\bf 0}\\
	-\dfrac{\delta t}{h\delta \rho}&\ddots &{\bf 0}\\
	{\bf 0}&\ddots& 1+\dfrac{\delta t}{h\delta \rho}
   \end{array}
   \right).
$

Using the  parameters $T=1$, $L=1$, $\alpha = 0.5$, $\delta t = 0.001$, $\delta x =0.01$, $\delta \rho = 0.001$
and initial conditions, $y_0(x)=1-\cos(2\pi x)$ and $z_0(\rho)=0.1\sin(-2\pi \rho h)$ we obtain the following figure, that represents $t\mapsto \ln(E(t))$ for different values of $h$ and $\beta$. We can see that when there is no delay, the energy is quickly exponentially decreasing and if the coefficient of delay, $\beta$, or the time delay, $h$ increase, then the energy is less exponentially decreasing.

\begin{center} \label{fig}
\includegraphics[scale=0.7]{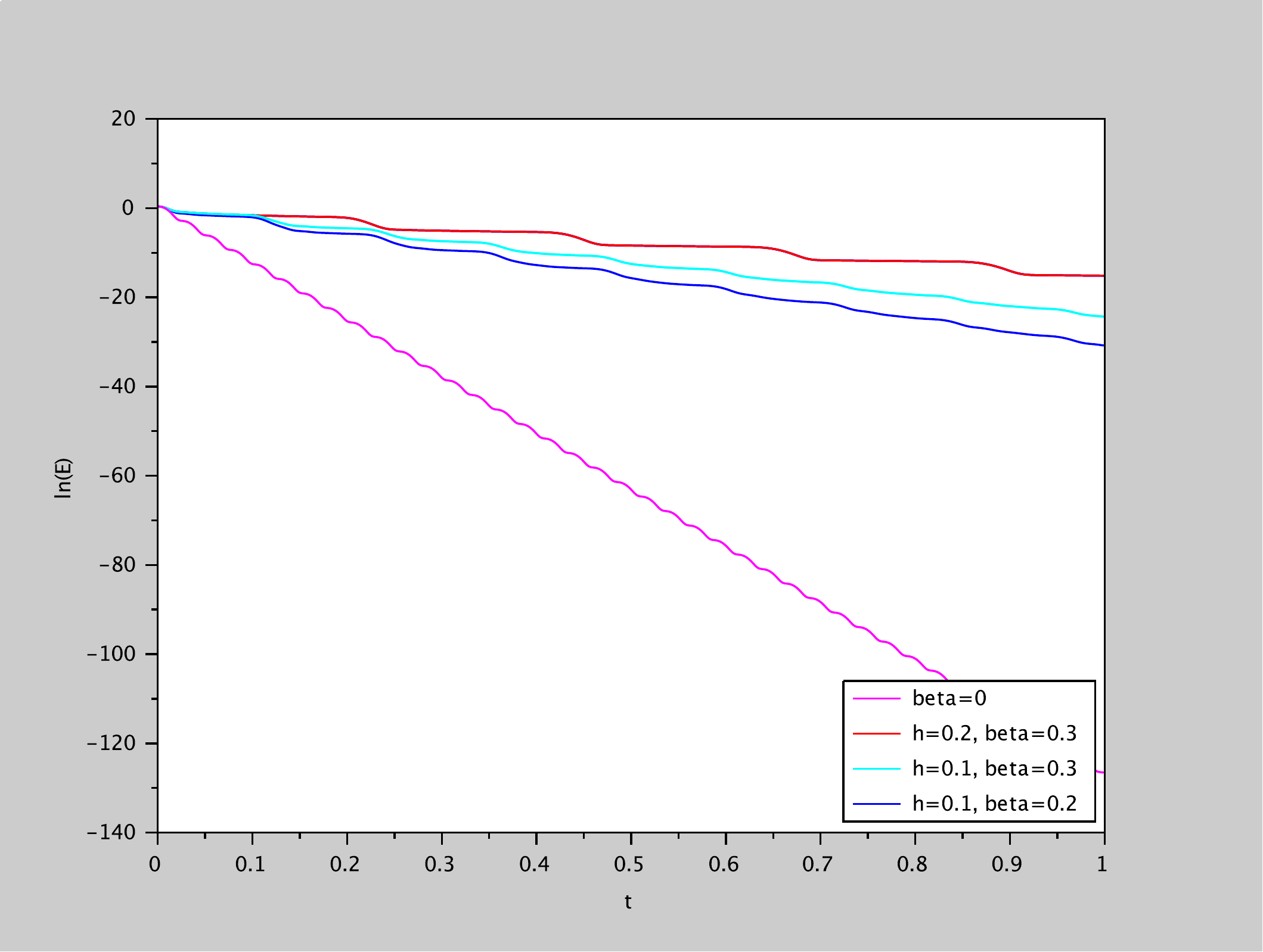}
\end{center}

\section{Conclusion}
In this article, we presented two different methodological approaches to prove some exponential stability results for the nonlinear KdV equation with a delayed term. More precisely, the KdV equation under study contains a boundary feedback, partly time delayed.
The main difficulty that our work faces is the nonlinearity of the KdV equation confronted with a delay term. The first main theorem proves a stability result with a quantified decay rate using a Lyapunov functional built to deal with both difficulties. It faces a technical limitation on the length of the spatial domain ($L< \pi\sqrt{3}$), more restrictive than the first critical length that KdV equations know.  The second stability result holds for any non-critical length, but without yielding any information on the exponential decay rate of the energy, since it follows a contradiction argument confronting an observability estimate.

A subject of future research could be the study of the stability of the nonlinear KdV equation specifically for critical lengths, using for instance a localised feedback as in \cite{Perla_2002}, but with some delay.



\end{document}